\documentclass[11pt]{article}

\usepackage{amssymb}
\usepackage{latexsym}
\newtheorem{theorem}{Theorem}[section]
\newtheorem{lem}[theorem]{Lemma}

\newtheorem{conj}[theorem]{Conjecture}

\newenvironment{proof}{\noindent {\it Proof.~~}\ }{\ \mbox{$\Box$}}

\newcommand{\subgp}[1]{\langle{#1}\rangle}
\newcommand{\beeq}{\begin{eqnarray*}}
\newcommand{\eneq}{\end{eqnarray*}}

\def\Z{\mathbb Z}

\begin{document}
\title{On complete subsets of the cyclic group}
\author{ {Y. O. Hamidoune}\thanks{
Universit\'e Pierre et Marie Curie, E. Combinatoire, Case 189, 4
Place Jussieu, 75005 Paris, France. \texttt{yha@ccr.jussieu.fr}}
\and {A.S. Llad\'o}\thanks{
 Universitat
Polit\`ecnica de Catalunya, Dept. Matem\`atica Apl. IV;  Jordi
Girona, 1, E-08034 Barcelona, Spain. \texttt{allado@ma4.upc.edu}}
\and {O. Serra}\thanks{
 Universitat
Polit\`ecnica de Catalunya, Dept. Matem\`atica Apl. IV;  Jordi
Girona, 1, E-08034 Barcelona, Spain. \texttt{oserra@ma4.upc.edu}} }

\date{}

\maketitle

\begin{abstract}
A subset $X$ of an abelian $G$ is said to be {\em complete} if every
element of the subgroup generated by $X$ can be expressed as a
nonempty sum of distinct elements  from $X$.

Let $A\subset \Z_n$ be such that all the elements of $A$ are coprime with $n$.
Solving a conjecture of
Erd\H{o}s and Heilbronn, Olson proved that
 $A$ is complete if $n$ is a prime and if  $|A|>2\sqrt{n}.$
 Recently  Vu proved that there is an absolute constant $c$,
 such that for an arbitrary large  $n$, $A$ is complete if $|A|\ge c\sqrt{n},$
 and conjectured that $2$ is  essentially the right value of $c$.

 We show that $A$ is complete if $|A|> 1+2\sqrt{n-4}$, thus proving
 the last conjecture.
\end{abstract}

\section{Introduction}

The additive group of integers modulo $n$ will be denoted by $\Z_n$.

Let $G$ be a finite Abelian group and let $X\subset G$. The subgroup
generated by  a subset $X$ of $G$ will be denoted $\subgp{ X}$. For
a positive integer $k$, we shall write $$k\wedge X =\left\{
\sum_{x\in A} x \ \ \ \Big| \ \ \ A\subset X \mbox{ and  }
|A|=k\right\}.$$ Following the terminology of \cite{tv} we write
$$
S_X=\bigcup _{k\ge 1} k\wedge X.
$$
The set  $X$  is said to be {\em complete} if $S_X=\subgp{X}.$ The
reader may find the connection between this notion and the
corresponding notion for integers in \cite{tv}. We shall also write
$$S_X^0=S_X\cup \{ 0\}.$$ Note that $S_X^0=\sum_{x\in X} \{ 0,x\}$.

Let $p$ denote a prime number and let $A\subset \Z_p\setminus \{
0\}$. Erd\H{o}s and Heilbronn \cite{EH} showed that  $A$ is complete
if $|A|\ge \sqrt{18}\sqrt{p}$, and conjectured that $\sqrt{18}$ can
be replaced by $2$. This conjecture was proved by Olson\cite{OLS1}.
More precisely, Olson's Theorem states that $A$ is complete if
$|A|\ge \sqrt{4p-4}.$ This result was sharpened by Dias da Silva and
one of the authors \cite{DH} by showing that $|k\wedge A|=p$,  if
$|A|\ge \sqrt{4p-4},$ where  $k=\lceil \sqrt{p-1}\ \rceil$. They
also showed that $|(j\wedge A)\cup ((j+1)\wedge A)|=p$,   if $|A|\ge
\sqrt{4p-8},$ where $j=\lceil \sqrt{p-2}\ \rceil$.

Let $G$ be a finite abelian group and let $A\subset G\setminus \{
0\}.$ Complete sets  for  general abelian group  were investigated
by Diderrich and Mann \cite{DM}. Diderrich  \cite{DI} proved that,
if $|G|=pq$ is the product of two primes, then $A$ is complete if
$|A|\ge p+q-1.$

Let $p$ be the smallest prime dividing $|G|.$ Diderrich conjectured
\cite{DI} that $A$ is complete, if $|G|/p$ is composite and
$|A|=p+|G|/p-2.$ This conjecture was  finally proved by Gao and one
of the authors \cite{GH}. More precise results were later proved by
Gao and the present authors \cite{GAOY}.
Note that the bound of Diderrich is best possible, since one may construct
 non complete sets of size $p+|G|/p-3$.

 However the result of Olson was extended recently by Vu
\cite{vu} to general cyclic groups.
Let $A\subset \Z_n$ be such that all the elements of $A$ are coprime with $n$.
 Vu proved that there is an absolute constant $c$ such that, for an arbitrary large
 $n$, $A$ is complete if $|A|\ge c\sqrt{n}.$
 The proof of Vu is rather short { and depends on a recent result of
 Szemer\'edi and  Vu \cite{sv}.} In the same paper
 Vu conjectures that the constant is essentially  $2$.

Our main result is the following:

\begin{theorem} \label{mainth}
Let $A$ be a subset of $\Z_n$ be such that all the elements of $A$
are coprime with $n$. If $|A|>1+2\sqrt{n-4}$ then $A$ is complete.

\end{theorem}

This result implies the validity of the last conjecture of Vu. We
conjecture the following:

\begin{conj} \label{mainconj}
Let $A\subset \Z_n$ be such that all the elements of $A$ are coprime
with $n$ and $|A|\ge \sqrt{4n-4}$. Then $|k\wedge A|=n$, where
$k=\lceil \sqrt{n-1}\ \rceil$.
\end{conj}

\section{Some tools}

In this section we present known material and some easy applications
of it. We give short proofs in order to make the paper
self-contained.

Recall the following well-known and easy lemma.

\begin{lem}  Let $G$ be a  finite group.
Let $X$ and $Y$ be subsets of $G$ such that
 $X+Y\neq G.$  Then $|X|+|Y|\leq |G|$.
\label{prehistorical}
\end{lem}

\begin{proof}
Take $a\in G\setminus (X+Y)$. We have $(a-Y)\cap X=\emptyset$.
\end{proof}

We use also the Chowla's Theorem \cite{MAN,NAT} :

\begin{theorem}[Chowla \cite{MAN,NAT}]
 Let $n$ be a positive integer and
let $X$ and $Y$ be non-empty subsets of $\Z_n$. Assume that $0\in Y$
and that the elements of $Y\setminus \{ 0\}$ are coprime with $n$.
 Then $$|X+Y|\geq \min (n,|X|+|Y|-1).$$\label{CD}
\end{theorem}

\begin{proof}
The proof is by induction on $|Y|$, the result being obvious for
$|Y|=1.$ Assume first that $Y \subset X-x,$ for all $x\in X.$ Then
$X+Y\subset X$, and hence $X+Y= X$. It follows that $X+Y=X+nY=\Z_n.$

Assume now that $Y\not\subset X-x,$ for some $x\in X.$ Then $0\in
Y\cap (X-x)$ and  $|Y\cap (X-x)|<|Y|.$ By the induction hypothesis,
$|X|+|Y|-1\le |((X-x)\cup Y)+((X-x)\cap Y)|\le |(X-x)+Y|.$
\end{proof}

Let $B\subset G$ and  $x\in G$. Following Olson, we write
$$
\lambda_B(x) = |(B+x)\setminus {B}|.
$$

The following result is implicit in \cite{OLS1}:
\begin{lem}[Olson, \cite{OLS1}] \label{implicit}
Let $Y$ be a nonempty subset of $G\setminus \{ 0\}$, $z\notin Y$ and
$y\in Y$. Put  $B=S_Y^0$. Then
\begin{equation}
|B| \geq |S_{Y\setminus \{ y\}}^0|+ \lambda_B(y), \label{eqolson1}
\end{equation}
and
\begin{equation}
|S_{Y\cup \{z\} }^0| = |S_{Y}^0|+ \lambda_B(z). \label{eqolson2}
\end{equation}
\end{lem}

\begin{proof}
Clearly we have $B\cap (\overline{B}-y)\subset B\setminus
S_{Y\setminus \{ y\}}^0$ and hence $\lambda_B (y)=|B\cap
(\overline{B}-y)|\le | B|-| S_{Y\setminus \{ y\}}^0|$.

>From $S_{Y\cup \{ z\}}^0=B+\{ 0,z\}$ we have $|S_{Y\cup \{
z\}}^0|=|B|+|(B+z)\setminus B|\}=|B|+\lambda _B(z).$
\end{proof}

We need the following helpful result also due  to Olson:

\begin{lem}[Olson \cite{OLS1}] Let $B$ and $C$ be nonempty subsets of
an abelian group $G$ such that $0\not\in C$. Then,
  \begin{eqnarray}
   \lambda_B(x) &=& \lambda_B(-x).  \label{eq:-x}\\
 \lambda_B(x+y) &\leq& \lambda_B(x)+\lambda_B(y).\label{eq:x+y}\\
 \sum _{x\in C} \lambda_B(x)&\geq& |B|(|C|-|B|+1). \label{eq:clique}
  \end{eqnarray}
\end{lem}

\begin{proof}
For each $x\in G$ we have
\begin{eqnarray*}
 |(B+x)\cap \overline{B}|&=&|B+x|-|(B+x)\cap B|\\
 &=&|B-x|-|B\cap (B-x)|\\
&=&|\overline{B}\cap (B-x)|=\lambda_B(-x),
\end{eqnarray*}
proving (\ref{eq:-x}). Let $x,y\in G$. Then,
\begin{eqnarray*}
\lambda_B(x+y) &=& |(B+x+y)\cap \overline{B}|\\
&=& |(B+x)\cap (\overline{B}-y)|\\
&=& |(B+x)\cap \overline{B}\cap (\overline{B}-y)|+
|(B+x)\cap B \cap (\overline{B}-y)|\\
&\le& |(B+x)\cap \overline{B}|+
| B\cap (\overline{B}-y)|\\
&=& \lambda_B(x)+\lambda_B(y),
\end{eqnarray*}
proving (\ref{eq:x+y}). Finally,
\begin{eqnarray*}
  \sum _{x\in C} \lambda_B(x)&\geq& \sum _{x\in C} (|B+x|-|B\cap (B+x)|)\\
 &\geq& |C| |B|-\sum _{x\in C}|B\cap (B+x)|\\
 &\geq& |C| |B|-\sum _{x\in G\setminus 0}|B\cap (B+x)|\\
 &=&|B|(|C|-|B|+1) ,
  \end{eqnarray*}
  proving (\ref{eq:clique}).
 \end{proof}

 \section{The main result}

The next Lemma is the key tool for our main result.

\begin{lem}\label{lem:eh} Let $A$ and $B$ be nonempty subsets of
$\Z_n$. Assume that $A\cap (-A)=\emptyset$ and that each element in
$A$ is coprime with $n$. Put $a=|A|$ and $b=|B|$. Assume also that
$a\ge 3$ and  $2b\le n+2$. Then
\begin{equation}\label{eq:lamb}
 \max_{x\in A}\lambda_B(x)>  a- \frac{a(a-3)}{b}.
\end{equation}
In particular, if $2b\ge a(a-3)$, then
\begin{equation} \label{eq:plus}
  \max_{x\in A}\lambda_B(x)\ge a-1.
\end{equation}
\end{lem}

\begin{proof} Put $A^*=A\cup(-A)\cup\{ 0\}$. Let $t<n$ be a positive integer
and set $$t=2ma+r,\; m\ge 0,\;  0\le r\le 2a-1.$$
 Let
$C_j=jA^*$. By Chowla's theorem, $|C_j|\ge \min \{ n,
2ja+1\}=2ja+1$, for $j\le m$. Therefore we can choose a set
$C\supset A^*$  of cardinality $t+1$ which intersects $C_j$ in
exactly  $2ja$ elements $j=2,\ldots ,m$, and intersects $C_{m+1}$ in
exactly $r$ elements. Let $E=C\setminus \{ 0\}$. Let
$\alpha=\max\{\lambda_B(x):\; x\in A\}$. By (\ref{eq:-x}) we have
$\lambda_B(x)\le \alpha ,$ for all $x\in A^*$. For an element $x$ in
 $C_j$ there are elements $x_1, \cdots , x_j\in A^*$  such that
$x=x_1+\cdots +x_j$. In view of  (\ref{eq:x+y}) we have
$\lambda_B(x)\le \lambda (x_1)+ \cdots +\lambda (x_j)\le j \alpha .$
Therefore,
\begin{eqnarray*}\sum_{x\in E}\lambda_B(x)&\le& \alpha 2a+2\alpha 2a+\cdots +m\alpha 2a+r(m+1)\alpha\\
&=& \alpha(m+1)(ma+r)=\frac{\alpha(t-r+2a)(t+r)}{4a}\\
&\le &\frac{\alpha(t+a)^2}{4a}.\end{eqnarray*} By using
(\ref{eq:clique}) we have
$$
\alpha\ge \frac{4a\sum_{x\in E}\lambda_B(x)}{(t+a)^2}
\ge \frac{4ab(t-b+1)}{(t+a)^2}.
$$
In particular, since   $2b\le n+2$, we can set  $t=2b-3$ to get,
\begin{eqnarray*}
\alpha&\ge&  \frac{4ab(b-2)}{(2b+a-3)^2}\\
&\ge&  \frac{a(b-2)}{b}(1-\frac{a-3}{b})\\
&>& a- \frac{a(a-3)}{b},
\end{eqnarray*}
where we have used $a\ge 3.$  In particular, if $2b\ge a(a-3)$, then
$\alpha>a-2$ so that $\alpha\ge a-1$. This completes the proof.
\end{proof}

Lemma \ref{lem:eh} gives the following estimation for the
cardinality of the set of subset sums.

\begin{lem} Let  $A\subset \Z_n$ such that $A\cap (-A)=\emptyset$ and
every element of $A$ is coprime with $n$. Also assume $|A|\ge 2.$ Then
$$|S_A^0|\ge \min \{ \frac{n+2}{2}, 3+\frac{|A|(|A|-1)}{2}\}.$$
\label{mainlemma}
\end{lem}

\begin{proof} We shall prove the result by induction on $a=|A|$, the result being obvious for
$a=2$. Suppose $a>2$. Put $B=S_A^0$. We may assume $b=|B|\le
\frac{n}{2}+1$ so that $2b\le n+2$. By the induction hypothesis,
$2b\ge 6+(a-1)(a-2)>a(a-3)$.

By (\ref{eq:plus}) there is an $x\in A$ with $\lambda_B(x)\ge a-1$.
Then, by Lemma \ref{implicit},
$$
|B|\ge |S_{A\setminus \{ x\}}^0|+\lambda_B(x)\ge 3+(a-2)(a-1)/2+a-1=
3+\frac{a(a-1)}{2},$$ as claimed.\end{proof}

We are now ready for the proof of Theorem \ref{mainth}.


 {\it Proof
of Theorem \ref{mainth}.\/} Suppose $A$ non complete and put
$|A|=k$.
 Let $X,Y$ be disjoint subsets of $A$. We
clearly have $S_{X}+S_{Y}^0\subset S_A\neq  \Z_n$.  Since
$|S_{X}|\ge |S_{X}^0|-1,$ we have
\begin{equation}\label{eqfinal}
|S_{X}^0|+|S_{Y}^0|\leq n+1, \end{equation}
 by   Lemma
\ref{prehistorical}.

Partition $A=A_1\cup A_2$ into two almost equal parts, i.e.
$|A_1|=\lceil k/2\rceil $ and $|A_2|=\lfloor k/2\rfloor $, such that
$A_i\cap (-A_i)=\emptyset$, $i=1,2$.

We must have \begin{equation}
 3+\lfloor\frac{k}{2}\rfloor(\lfloor\frac{k}{2}\rfloor-1)/2<
(n+2)/2,\label{eq:k}\end{equation} since otherwise, by Lemma
\ref{mainlemma}, we have  $|S_{A_1}^0|+|S_{A_2}^0|\ge n+2,$
contradicting (\ref{eqfinal}).

{\it Case 1.\/} $k$ even.

Then we have by (\ref{eq:k})
$$n/2>2+\frac{k}{2}(\frac{k}{2}-1)/2=2+k(k-2)/8,$$
and hence $(k-1)^2+16\le 4n,$ a contradiction.

{\it Case 2.\/} $k$ odd.

Put $a=\frac{k-1}{2}$.  In view of  (\ref{eq:k}), Lemma
\ref{mainlemma} implies $$|S_{A_2}^0|\ge 3+a(a-1)/2.$$ By
(\ref{eq:plus}) with $B=S_{A_2}^0$, there is a $y\in A_1$ such that
$$
\lambda_B(y)\ge a-1.
$$
Put $C_1=A_1\setminus \{y\}$ and $C_2=A_2\cup \{y\}$. Then we have,
by Lemma \ref{implicit},
$$
|S_{C_2}^0|\ge |S_{A_2}^0|+\lambda_B(y)\ge 3+a(a-1)/2+a-1=
2+\frac{a(a+1)}{2}.
$$
On the other hand, from (\ref{eq:k}) and   Lemma \ref{mainlemma} we
get
$$
|S_{C_1}^0|\ge 3+\frac{a(a-1)}{2}.
$$
By (\ref{eqfinal}), $$n+1\ge |S_{C_1}^0|+|S_{C_1}^0|\geq
3+a(a-1)/2+2+a(a+1)/2=5+a^2.$$ Therefore $4n\ge 16+4a^2=16+(k-1)^2$,
a contradiction. This completes the proof. $\Box$


\begin{thebibliography}{99}


\bibitem{DH} J.A. Dias da Silva and Y. O. Hamidoune,
Cyclic spaces for Grassmann derivatives and additive theory, {\it
      Bull. London Math. Soc.}, 26 (1994), 140-146.


\bibitem{DI} G.T. Diderrich, An addition theorem for abelian groups
of order pq, {\it  J. Number Theory} {\bf 7} (1975), 33-48.

\bibitem{DM} G. T. Diderrich and H. B. Mann, {\it Combinatorial
problems in finite abelian groups}, In: "A survey of Combinatorial
Theory" (J.L. Srivasta et al. Eds.), pp. 95- 100, North- Holland,
Amsterdam (1973).




\bibitem{EH} P. Erd\H{o}s and H. Heilbronn,
On the Addition of residue classes mod  $p$, {\it Acta Arith.} {\bf
9} (1964), 149-159.

\bibitem{GH} W. Gao and Y.O. Hamidoune, On  additive bases,
{\it Acta Arith.}  {\bf 88}  (1999), 3, 233-237.



\bibitem{GAOY} W. Gao, Y.O. Hamidoune
A. S. Llad\'o and O. Serra, Covering a finite abelian group by
subset sums.  {\it Combinatorica } {\bf  23 } (2003),  no. 4,
599--611.


\bibitem{MAN} H.B. Mann, {\it Addition Theorems},   R.E.
Krieger, New York, 1976.


\bibitem{OLS1} J. E. Olson,  An addition theorem mod $p$, {\it J. Comb.
Theory} {\bf  5} (1968), 45-52.




\bibitem{CHO} S.~Chowla, A theorem on the addition of residue
classes: applications to the number $\Gamma (k)$ in Waring's
problem, {\it Proc.Indian Acad. Sc.} {\bf 2} (1935) 242--243.

\bibitem{NAT}M. B. Nathanson,
{\it Additive Number Theory. Inverse problems and the geometry of
sumsets}, Grad. Texts in Math. 165, Springer, 1996.

\bibitem{sv} E. Szemer\'edi and V.H. Vu, Long arithmetic
progressions in finite and infinite sets, {\it Annals of Math.}, to
appear.

\bibitem{tv} T. Tao and V.H. Vu, {\it Additive Combinatorics}, Cambridge Studies in Advanced Mathematics 105
(2006), Cambridge Press University.

\bibitem{vu} V.H. Vu, Olson Theorem for cyclic groups, Preprint, arXiv:math.NT/0506483 v1, 23 june 2005.
\end{thebibliography}
\end{document}